\documentclass[fleqn]{mat01}
\usepackage{times,mathtimy,amssymb,epsfig}
\begin{document}

\setcounter{page}{171}
\firstpage{171}

\font\xx=msam5 at 10pt
\def\qed{\mbox{\xx{\char'03\!}}}

\makeatletter
\def\artpath#1{\def\@artpath{#1}}
\makeatother \artpath{C:/mathsci-arxiv/MAY-2003}

\font\xx=msam5 at 9pt
\def\ab{\mbox{\xx{\char'03}}}

\font\sa=tibi at 10.4pt


\def\defi{\trivlist\item[\hskip\labelsep{\bf DEFINITION.}]}
\def\remarks{\trivlist\item[\hskip\labelsep{\it Remarks}]}
\def\noot{\trivlist\item[\hskip\labelsep{{\it Note.}}]}
\def\thoe{\trivlist\item[\hskip\labelsep{{\bf Theorem}}]}
\newtheorem{theo}{Theorem}
\renewcommand\thetheo{\arabic{section}.\arabic{theo}}
\newtheorem{theor}[theo]{\bf Theorem}
\newtheorem{lem}[theo]{Lemma}
\newtheorem{propo}{\rm PROPOSITION}
\newtheorem{rema}[theo]{Remark}
\newtheorem{exam}{Example}
\newtheorem{coro}{\rm COROLLARY}

\newcommand{\Z}{\ensuremath{\mathbb Z}}
\newcommand{\C}{\ensuremath{\mathbb C}}
\newcommand{\R}{\ensuremath{\mathbb R}}
\newcommand{\T}{\ensuremath{\mathbb T}}
\newcommand{\N}{\ensuremath{\mathbb N}}
\newcommand{\G}{\ensuremath{\mathcal G}}
\newcommand{\M}{\ensuremath{\mathcal M}}
\renewcommand{\L}{\ensuremath{\mathcal L}}





\title{Wavelet subspaces invariant under groups of translation operators}

\markboth{Biswaranjan Behera and Shobha Madan}{Translation invariant wavelet subspaces}

\author{BISWARANJAN BEHERA$^{*}$ and SHOBHA MADAN}

\address{Department of Mathematics, Indian Institute of Technology, Kanpur 208 016, India\\
\noindent $^{*}${\it Current address}: Statistics and Mathematics Unit, Indian Statistical Institute,\\
\noindent 203 B.T. Road, Kolkata 700 108, India\\
\noindent E-mail: br\_behera@yahoo.com; madan@iitk.ac.in.}

\volume{113}

\mon{May}

\parts{2}

\Date{MS received 6 November 2002}

\begin{abstract}
We study the action of translation operators on wavelet
subspaces. This action gives rise to an equivalence relation on the set
of all wavelets. We show by explicit construction that each of the
associated equivalence classes is non-empty.
\end{abstract}

\keyword{Wavelet; multiresolution analysis; multiresolution analysis
wavelet; minimally supported frequency wavelet; wave\-let set;
translation invariance.}

\maketitle

\section{Introduction}

A function $\psi\in L^{2}(\R)$ is said to be a {\it wavelet} if the
system of functions $\{\psi_{j,k}:j,k\in \Z \}$ forms an orthonormal
basis for $L^{2}(\R)$, where
\begin{equation*}
\psi_{j,k}(x)=2^{j/2}\psi(2^jx-k),\quad  j,k\in \Z.
\end{equation*}
Mallat~\cite{Mal} and Meyer~\cite{Mey} provided a framework to construct
wavelets  through the concept of multiresolution analysis (MRA). A
sequence $\{V_j:j\in\Z\}$ of closed subspaces of $L^2(\R)$ is called an
MRA if the following conditions hold:

\begin{enumerate}
\leftskip .4pc
\renewcommand{\labelenumi}{(\roman{enumi})}
\item $V_{j}\subset V_{j+1}$ for all $j\in\Z$,
\item $f\in V_{j}$ if and only if $f(2\cdot)\in V_{j+1}$ for all
$j\in\Z$,
\item $\bigcup_{j\in\Z} V_{j}$ is dense in $L^{2}(\R)$,
\item $\bigcap_{j\in\Z} V_j = \{0\}$, and
\item there exists a function $\varphi\in V_{0}$, called the {\it scaling
function}, such that $\{\varphi(\cdot-k):k\in\Z\}$ forms an orthonormal
basis for $V_{0}$.
\end{enumerate}

One can always construct a wavelet from an MRA (see, for instance,
Ch.~2 in~\cite{hw}), but not every wavelet can be obtained in this
manner. The first example of a wavelet which cannot be obtained from an
MRA was given by Journ\'e. In~\cite{BM} we characterized
a large class of wavelets, which also includes Journ\'e's wavelet, and
proved that none of them is associated with an MRA (see
Theorem~\ref{thm:sn}).

The following theorem, which characterizes all wavelets of $L^{2}(\R)$,
was proved independently by Gripenberg~\cite{grip} and
Wang~\cite{wan} (see also~\cite{hkls,hw}).

\begin{theor}[\!]
\label{thm:wavelet1}
Let $\psi\in L^2(\R)$ with $\|\psi\|_2=1$. Then $\psi$ is a wavelet of
$L^2(\R)$ if and only if
\begin{align}
 & \sum_{j\in \Z}|\hat\psi(2^j\xi)|^2=1\quad\mbox{for a.e.}~\xi\in\R.
   \label{eqn:1W1} \\
 & \sum_{j\geq 0} \hat\psi(2^j\xi)\overline{\hat\psi(2^j(\xi +2m\pi))}=0
   \quad\mbox{for a.e.}~\xi\in\R,~\mbox{for all}~m\in 2\Z+1.
   \label{eqn:1W2}
\end{align}
\end{theor}

We use the following definition of the Fourier transform:
\begin{equation*}
\hat f(\xi)=\int_{\R}f(x)\hbox{e}^{-i\xi x}\hbox{d}x,\quad\xi\in\R.
\end{equation*}

If $\psi$ is a wavelet, then the support of $\hat\psi$ must have measure
at least $2\pi$. This minimal measure is achieved if and only if
$|\hat\psi|$ is the characteristic function of some measurable subset
$K$ of $\R$. Such a wavelet is called a minimally supported frequency
(MSF) wavelet and the associated set $K$ is called a {\it wavelet set}.
We refer to~\cite{hw} for proofs of the above statements.

As we mentioned earlier, a wavelet need not be associated with an MRA,
but it can be made to be associated with an MRA-like structure in the
following manner. Given a wavelet $\psi$ we define the closed subspaces
$V_j$, $j\in\Z$, by
\begin{equation}
\label{E.vjs}
V_{j} = \overline{{\rm span}}\{\psi_{l,k}:l<j, k\in\Z\}, \quad j\in\Z.
\end{equation}
It is easy to verify that these subspaces satisfy properties (i)--(iv)
of an MRA. Moreover, instead of (v), the following weaker property
holds:\vspace{.3pc}

\begin{enumerate}
\leftskip .2pc
\item[(v$'$)] $V_{0}$ is invariant under translation by integers.
\end{enumerate}

A sequence of closed subspaces $\{V_{j}:j\in\Z\}$ of $L^{2}(\R)$ which
satisfies properties (i)--(iv) and ($\hbox{v}'$) is called a generalized
MRA (GMRA).

Madych (\S3 in~\cite{Mad}) characterized all MRAs for which each
$V_{j}$ is invariant under translation by all real numbers. Let
$T_\alpha$, $\alpha\in\R$, be the translation operator defined by
$T_\alpha f(x)=f(x-\alpha)$. An MRA $\{V_{j}:j\in\Z\}$ is called a {\it
translation invariant} MRA if $T_\alpha(V_j)\subset V_j$ for all
$j\in\Z$ and for all $\alpha\in\R$. Madych proved that the only
translation invariant MRAs are those for which the Fourier transform of
the scaling function is the characteristic function of a set. In other
words, the associated wavelet is an MSF wavelet.

Walter~\cite{Wal1,Wal2} modified the definition of translation
invariance to include wavelets other than the MSF wavelets. An MRA is
called {\it weakly translation invariant} if $T_\alpha(V_j)\subset
V_{j+1}$ for $j\in\Z$ and $\alpha\in\R$. Clearly, every translation
invariant MRA is weakly translation invariant, since $V_j\subset
V_{j+1}$. Walter gave necessary conditions for an MRA to be weakly
translation invariant. With an additional condition, these conditions
also turned out to be sufficient for the weak translation invariance of
MRAs associated with a class of Meyer-type wavelets (see \S10.5
in~\cite{Wal2}). Thus, we see that the notion of translation invariance
and weak translation invariance are applicable only to a limited class
of wavelets.

In this article we investigate the translation invariance from a
slightly different point of view. We work in the more general set up of
GMRAs in order to include {\it all} wavelets. We still demand the spaces
$V_j$ to be invariant under translations, but only by dyadic rationals
at a fixed level instead of translation by all reals. More precisely, we
ask the following question:

Let $n\in\N$. Does there exist a wavelet such that the space
$V_{0}$ of the associated GMRA is invariant under translation by elements
of the form $m/2^{n}$ for all $m\in\Z$?

The purpose of this paper is to give an affirmative answer to this
question. We do this by explicitly constructing such wavelets for each
$n$.

Let $\G_{n}$, $n\in\N\cup\{0\}$, and $\G_\infty$ be the following groups
of (unitary) translation operators:
\begin{equation*}
\G_n=\{T_{{m}/{2^n}}:m\in\Z\},\quad
\G_\infty=\{T_\alpha:\alpha\in\R\}.
\end{equation*}
Denote by $\L_{n}$ the collection of all wavelets $\psi$ such that the
space $V_{0}$ (of the GMRA $\{V_{j}\}$ associated with $\psi$) is invariant
under the group $\G_n$. Clearly, $\L_0$ is the set of all wavelets, and
\begin{equation*}
 \L_{0}\supset\L_{1}\supset\L_{2}\supset\cdots\supset\L_{n}
 \supset\L_{n+1}\supset\cdots\supset\L_{\infty}.
\end{equation*}
These inclusions naturally defines an equivalence relation on the set of
all wavelets. The equivalence classes are given by \hbox{$\M_n = \L_n\!\setminus\! \,\L_{n+1}$}, $n\in\N\cup\{0\}$, and $\M_\infty = \L_\infty$. Therefore,
$\M_n$, $n\in\N\cup\{0\}$, is the class of wavelets such that $V_0$ is
invariant under the group $\G_n$ but not under $\G_{n+1}$. This
equivalence relation was defined by Weber in~\cite{web} and
the equivalence classes were characterized in terms of the support of
the Fourier transform of the wavelets.

Given a wavelet $\psi$, let $E(\psi,k)={\rm supp}~\hat\psi\cap({\rm
supp}~\hat\psi+2k\pi)$, $k\in\Z$, and ${\mathcal
E}(\psi)=\{k\in\Z:E(\psi,k)\not=\emptyset\}$. In this notation, the
characterization of $\M_n$, $n\in\N\cup\{0,\infty\}$, as given
in~\cite{web}, is the following.

\begin{theor}[\!]
\label{T.char}
{\rm (a)} $\M_\infty$ is precisely the collection of all MSF wavelets. {\rm (b)} The
equivalence class $\M_n,$ $n\in\N,$ consists of all wavelets $\psi$ such
that every element of ${\mathcal E}(\psi)$ is divisible by $2^n$ but
there exists an element of ${\mathcal E}(\psi)$ not divisible by
$2^{n+1}$. {\rm (c)} A wavelet $\psi$ belongs to $\M_0$ if and only if
${\mathcal E}(\psi)$ contains an odd integer.
\end{theor}

In the same paper, Weber produced examples of wavelets belonging to the
first few equivalence classes $\M_n$, namely for $n=0,1,2$ and $3$. This
motivated us to construct wavelets for each $\M_n$. After we constructed
these wavelets we came to know about the article~\cite{sw}, in which
Schaffer and Weber also constructed wavelets in $\M_n,n\in\N$,
using the method of `operator interpolation' (see chapter 5 in~\cite{dl}).

In this paper we will construct wavelets belonging to each of these
classes by a different method. Our approach is simpler than that
of~\cite{sw} in the sense that for each integer $n\geq 3$, we construct
a function $\psi_n$ such that $\psi_n$ has the required properties to be
in $\M_{n-2}$ as characterized in Theorem~\ref{T.char}. Then we show
that $\psi_n$ is a wavelet. This is the content of \S2. In
\S3 we construct a family of wavelets belonging to the equivalence
class $\M_{0}$.

\font\rnsi=cmbsy10 at 10pt
\font\sfdf=mtmib at 10.4pt
\font\sfff=mtmib at 8pt
\def\caM{\mbox{\rnsi \char'115}}

\section{Construction of wavelets in \protect{$\caM_{n}$, $\hbox{\it n}\geq \hbox{\bf 1}$}}
\label{subsec:equiv1}

Let $n\geq 2$ be an integer. Put
\begin{equation*}
\begin{array}{ll}
a_n=\displaystyle\frac{2^{n-1}}{2^n-1}\pi,        & b_n=2a_n=\displaystyle\frac{2^n}{2^n-1}\pi, \\[.9pc]
c_n=\displaystyle\frac{2^{n-1}(2^{n}-2)}{2^{n}-1}\pi, & d_n=2^n a_n =\displaystyle\frac{2^{2n-1}}{2^n-1}\pi,\\[.9pc]
e_n=\displaystyle\frac{2^{n}-2}{2^n-1}\pi.          &
\end{array}
\end{equation*}
Define $S_n= S_n^{+}\cup S_n^{-}$, where
$S_n^{+}=[a_n,b_n]\cup[c_n,d_n]$ and $S_n^- = -(S_n^+)$.

\begin{figure}[b]\vspace{.6pc}
\centerline{\epsfxsize=9cm\epsfbox{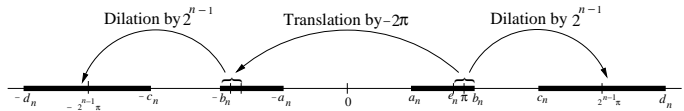}}\vspace{-.5pc}
\caption{The set $S_n$.}
\label{Fig:setsn}
\end{figure}

Observe that
\begin{equation}
\label{E.rel1}
[-b_n,-e_n]=[e_n,b_n]-2\pi,\ \quad [c_n,d_n]=2^{n-1}[e_n,b_n]
\end{equation}
and
\begin{equation}
\label{E.rel2}
[-d_n,-c_n]=2^{n-1}([e_n,b_n]-2\pi).
\end{equation}

The following theorem, proved by the authors in~\cite{BM}, characterizes
all wavelets $\psi$ of $L^2(\R)$ such that the support of $\hat\psi$ is
contained in the set $S_n$.

\setcounter{theo}{0}

\begin{theor}[\!]
\label{thm:sn}
Let $n\geq 2,$ $\psi\in L^2(\R),$ {\rm supp}~$\hat\psi\subseteq S_n$ and
$b(\xi)=|\hat\psi(\xi)|$. Then $\psi$ is a wavelet for $L^2(\R )$ if and
only if
\end{theor}\vspace{-.6pc}
{\it \begin{enumerate}
\leftskip .4pc
\renewcommand{\labelenumi}{{\rm (\roman{enumi})}}
\item $b(\xi)=1$ \quad for a.e. $\xi\in[a_n,e_n]\cup[-e_n,-a_n],$\vspace{.2pc}
\item $b^2(\xi)+b^2(2^{n-1}\xi)=1$ \quad for a.e. $\xi\in [e_n,b_n],$\vspace{.2pc}
\item $b^2(\xi)+b^2(\xi-2\pi)=1$ \quad for a.e. $\xi\in[e_n,b_n],$\vspace{.2pc}
\item  $b(\xi)=b(2^{n-1}(\xi-2\pi))$ \quad for a.e. $\xi\in[e_n,b_n],$\vspace{.2pc}
\item  $\hat \psi(\xi)={\rm e}^{i\theta(\xi)}b(\xi),$ where $\theta$ satisfies\vspace{-1pc}
\end{enumerate}
\begin{equation*}
\theta(\xi)+\theta(2^{n-1}(\xi-2\pi))-\theta(\xi-2\pi)-\theta(2^{n-1}\xi)
=(2m(\xi)+1)\pi,
\end{equation*}
for some $m(\xi)\in\Z,$ for a.e.
$\xi\in[e_n,b_n]\cap({\rm supp}~b)\cap(\frac{1}{2^{n-1}}{\rm supp}~b)$.
Moreover{\rm ,} if $n\ge 3${\rm ,} then none of these wavelets is associated with
an MRA.}

\begin{rema}
{\rm (a) It follows from Theorem~\ref{thm:sn} that $\hat\psi$ is completely
determined by its values on $[e_n,b_n]$. On this set, let $b=|\hat\psi|$
be an arbitrary
measurable function taking values between 0 and 1. Then, in view of
(\ref{E.rel1}) and (\ref{E.rel2}), $|\hat\psi|$ can be extented to other
sets of $S_n$ using properties (i)--(iv) (see Figure~\ref{Fig:setsn}).
Any function $\theta$ satisfying (v) now completely defines $\hat\psi$
on $\R$.
(b) The function $\theta(\xi)=2^{-(n-1)}\xi$ is a solution of the functional equation in Theorem~\ref{thm:sn}(v).
(c) If $({\rm supp}~b)\cap((1/2^{n-1}){\rm supp}~b)$ has an empty interior
in $[e_n,b_n]$, then $\theta$ can be chosen to be any measurable
function. In particular, we can take\break $\theta(\xi)=0$.}
\end{rema}

Let $n\geq 3$. By choosing $b=0$ a.e. on the interval $[e_n,b_n]$ and extending it to the other sets of $S_n$ we get a wavelet $\gamma_n$, where
\begin{equation*}
\hat\gamma_n=\chi_{_{W_n}}, \quad W_n=[-b_n,-a_n]\cup[a_n,e_n]\cup[c_n,d_n].
\end{equation*}

\begin{figure}[t]
\centerline{\epsfxsize=9cm\epsfbox{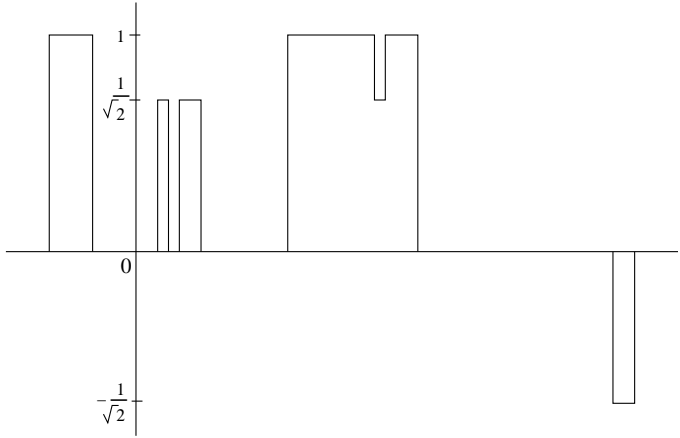}}\vspace{-.5pc}
\label{Fig:psihat} \caption{The function
$\hat\psi_n$.}\vspace{.6pc}
\end{figure}

We construct $\psi_n$ from this function in the following manner. We translate the interval
$[a_{n}/{2},e_{n}/{2}]+2^{n-1}\pi$ (which is a subset of $[c_n,d_n]$) to
the left by a factor of $2^{n-1}\pi$ and assign values $1/\sqrt{2}$ to
$\hat\psi_n$ on both these sets. Then, we translate $[a_n,e_n]$
to the right by a factor of $2^n\pi$ and assign to $\hat\psi_n$ the
value $1/\sqrt{2}$ on $[a_n,e_n]$ and $-1/\sqrt{2}$ on
$[a_n,e_n]+2^n\pi$. Assign the value $1$ to $\hat\psi_n$ on the
remaining sets  of
$W_n$ and $0$ elsewhere. More precisely, we have the following function
(see figure~\ref{Fig:psihat}):
\begin{equation*}
 \hat\psi_n(\xi)\!=\!
 \begin{cases}
 1 & {\rm if}~\xi\in[-b_n,-a_n]\cup[c_n,\frac{a_n}{2}+2^{n-1}\pi]
     \cup[\frac{e_n}{2}+2^{n-1}\pi,d_n],\\[.3pc]
 \frac{1}{\sqrt 2} & {\rm if}~\xi\in[\frac{a_n}{2},\frac{e_n}{2}]
     \cup[a_n,e_n]\cup[\frac{a_n}{2}+2^{n-1}\pi,
     \frac{e_n}{2}+2^{n-1}\pi],\\[.3pc]
 -\frac{1}{\sqrt 2} & {\rm if}~\xi\in[a_n+2^n\pi,e_n+2^n\pi],\\
 0 & {\rm otherwise.}
 \end{cases}
\end{equation*}$\left.\right.$\vspace{-1pc}

In the following lemma we list the translates and dilates of subsets of  supp~$\hat\psi_n$ which are again in supp~$\hat\psi_n$. The proof is a straightforward calculation and is  ommited.

\begin{lem}
\label{lemma2}
Let $F_n={\rm supp}~\hat\psi_n$.
\begin{enumerate}
\leftskip .8pc
\renewcommand{\labelenumi}{{\rm (\roman{enumi})}}
\item If $\xi\in[-b_n,-a_n],$ then $\xi+2k\pi\in F_n$ iff $k=0,$
and $2^j\xi\in F_n$ iff $j=0$.\vspace{.1pc}
\item If $\xi\in[a_{n}/{2},e_{n}/{2}],$
then $\xi+2k\pi\in F_n$ iff $k=0,2^{n-2},$ and $2^j\xi\in F_n$ iff $j=0,1$.\vspace{.1pc}
\item If $\xi\in[a_{n},e_{n}],$ then $\xi+2k\pi\in F_n$ iff $k=0,2^{n-1},$ and $2^j\xi\in F_n$ iff $j=0,-1$. \vspace{.1pc}
\item If $\xi\in[c_{n},(a_{n}/{2})+2^{n-1}\pi],$ then $\xi+2k\pi\in F_n$ iff $k=0,$ and $2^j\xi\in F_n$ iff $j=0$.\vspace{.1pc}
\item If $\xi\in[(a_{n}/{2})+2^{n-1}\pi,(e_{n}/{2})+2^{n-1}\pi],$
then $\xi+2k\pi\in F_n$ iff $k=0,-2^{n-2},$ and $2^j\xi\in F_n$ iff $j=0,1$.\vspace{.1pc}
\item If $\xi\in[(e_{n}/{2})+2^{n-1}\pi,d_{n}],$
then $\xi+2k\pi\in F_n$ iff $k=0,$ and $2^j\xi\in F_n$ iff $j=0$. \vspace{.1pc}
\item If $\xi\in[a_n+2^n\pi,e_n+2^n\pi],$
then $\xi+2k\pi\in F_n$ iff $k=0,-2^{n-1},$ and $2^j\xi\in F_n$ iff $j=0,-1$.
\end{enumerate}
\end{lem}

\begin{theor}[\!]
For each $n\geq 3,$ the function $\psi_n$ defined above is a wavelet and
belongs to the equivalence class $\M_{n-2}$.
\end{theor}

\begin{proof}
To prove that $\psi_n$ is a wavelet, it is sufficient to show that $\psi$ satisfies the following three properties (see Theorem~\ref{thm:wavelet1}):
\begin{enumerate}
\item[(a)] $\|\psi_n\|_2=1$, \item[(b)]
$\rho(\xi):=\sum_{j\in\Z}|\hat\psi_n(2^j\xi)|^2=1~ {\rm
for~a.e.}~\xi\in\R$, \item[(c)] $t_q(\xi):=\sum_{j\geq
0}\hat\psi_n(2^j\xi) \overline{\hat\psi_n(2^j(\xi
+2q\pi))}=0$~{\rm for~a.e.~}$\xi\in\R$, {\rm for~all~}$q\in
2\Z+1$.\vspace{-3pc}
\end{enumerate}
\end{proof}

\noindent{\it Proof of} (a).\hskip .2cm
We have
\begin{align*}
 \|\hat\psi_n\|_2^2
 &= (b_n-a_n)+\Bigl(\frac{a_n}{2}+2^{n-1}\pi-c_n\Bigr)+
       d_n-\Bigl(\frac{e_n}{2}+2^{n-1}\pi\Bigr) \\
 &\qquad +\frac{1}{2}\Bigl\{\left(\frac{e_n}{2}-\frac{a_n}{2}\right)
       +(e_n-a_n)+\Bigl(\frac{e_n}{2}-\frac{a_n}{2}\Bigr) +(e_n-a_n)\Bigr\} \\
 &= b_n-a_n+d_n-c_n+e_n-a_n=2\pi.
\end{align*}
Therefore,
$\|\psi_n\|_2=({1}/{2\pi})\|\hat\psi_n\|_2^2=1$.\vspace{.3cm}

\noindent{\it Proof of} (b).\hskip .2cm
Let $\xi>0$. Since $\rho(\xi)= \rho(2\xi)$, it is enough to show that
$\rho(\xi)=1$ for a.e.\ $\xi\in[\alpha,2\alpha]$ for some $\alpha >0$.
We will prove that $\rho(\xi)=1$ for a.e.\ $\xi\in
[a_n,2a_n]=[a_n,e_n]\cup[e_n,b_n]$.

Suppose $\xi\in [a_n,e_n]$. Then $2^j\xi\in F_n$ if and only if $j=0,-1$. So,
$\rho(\xi)=|\hat\psi_n(\xi)|^2+|\hat\psi_n({\xi}/{2})|^2
=({1}/{\sqrt2})^2+({1}/{\sqrt 2})^2=1$. Now,
$\xi\in [e_n,b_n]$ if and only if $2^{n-1}\xi\in [c_n,d_n]$. We write $[c_n,d_n]$ as a union of three intervals:
\begin{align*}
\hskip -1pc[c_n,d_n]
&=
 \left[c_n,\frac{a_n}{2}+2^{n-1}\pi\right]\!\cup\!
 \left[\frac{a_n}{2}+2^{n-1}\pi,\tfrac{e_n}{2}+2^{n-1}\pi\right]\!\cup\!
 \left[\frac{e_n}{2}+2^{n-1}\pi,d_n\right]\\
\hskip -1pc &=
 I_1\cup I_2\cup I_3, ~{\rm say}.
\end{align*}
If $2^{n-1}\xi\in (I_1\cup I_3)$, then $2^j(2^{n-1}\xi)\in F_n$ if and only if $j=0$
(see Lemma~\ref{lemma2}). So, $\rho(\xi)=|\hat\psi_n(2^{n-1}\xi)|^2=1$.
Also if $2^{n-1}\xi\in I_2$, then $2^j(2^{n-1}\xi)\in F_n$ if and only if $j=0$ or $1$, hence, $\rho(\xi)=|\hat\psi_n(2^{n-1}\xi)|^2+|\hat\psi_n(2^n\xi)|^2
=({1}/{\sqrt 2})^2+(-{1}/{\sqrt 2})^2=1$. Therefore, we get $\rho(\xi)=1$
for a.e.~$\xi>0$.

For $\xi<0,$ it suffices to show that $\rho(\xi)=1$ on $[-b_n,-a_n]$. On
this set, $2^j\xi\in F_n$ if and only if $j=0$. Hence,
$\rho(\xi)=|\hat\psi_n(\xi)|^2=1$ for a.e.\ $\xi\in[-b_n,-
a_n]$.\vspace{.3cm}

\noindent{\it Proof of} (c).\hskip .2cm
Since $t_q(\xi)=\overline{t_{-q}(\xi+2q\pi)}$, it is enough to show that
$t_q=0$ a.e., if $q$ is a negative odd integer. Suppose $q\not=-1,$ and
is odd. We have $2^jq\neq 0,\pm 2^{n-1},\pm 2^{n-2}$. Therefore, if
$2^j\xi\in F_n$, then by Lemma~\ref{lemma2} we observe that
$2^j\xi+2\cdot 2^jq\pi\not\in F_n$, which shows that each term of the
sum $t_q(\xi)$ is $0$. Hence, $t_q=0$ a.e.

It remains to prove that $t_{-1}(\xi)=0$ for a.e. $\xi\in\R$. We have
\begin{equation*}
t_{-1}(\xi)=\sum\limits_{j\geq 0}\hat\psi_n(2^j\xi)
\overline{\hat \psi_n (2^j\xi -2\cdot 2^j\pi)}.
\end{equation*}
By Lemma~\ref{lemma2} (see (vii) and (v)), we observe that both $2^j\xi$
and $2^j\xi-2\cdot 2^j\pi$ belong to $F_n$ only in the following two
cases:
\begin{enumerate}
\renewcommand{\labelenumi}{(\roman{enumi})}
\item $2^j\xi\in[a_n+2^n\pi, e_n+2^n\pi]$ and $j=n-1$.
\item $2^j\xi\in[\frac{a_n}{2}+2^{n-1}\pi, \frac{e_n}{2}+2^{n-1}\pi]$
and $j=n-2$.
\end{enumerate}
But both are equivalent to saying that $2^{n-1}\xi\in[a_n+2^n\pi,e_n+2^n\pi]$. So we get
$t_{-1}(\xi)=0$, if $2^{n-1}\xi\not\in[a_n+2^n\pi,e_n+2^n\pi]$.

Now, if  $2^{n-1}\xi\in[a_n+2^n\pi,e_n+2^n\pi]$, then
\begin{align*}
 t_{-1}(\xi) & =
 \hat\psi_n(2^{n-2}\xi)\overline{\hat\psi_n(2^{n-2}\xi-2^{n-1}\pi)}
 +\hat\psi_n(2^{n-1}\xi)\overline{\hat\psi_n(2^{n-1}\xi-2^n\pi)}\\
 & = \frac{1}{\sqrt2}\cdot \frac{1}{\sqrt2}+\left(-\frac{1}{\sqrt2}\right)
 \cdot\frac{1}{\sqrt2} \\
 & = 0.
\end{align*}
This completes the proof of (c). Therefore, $\psi_n$ is a wavelet.

Our claim now is that $\psi_n\in\M_{n-2}$. By referring to
Lemma~\ref{lemma2} again we observe that ${\mathcal E}(\psi_n)=\{0,\pm
2^{n-2},\pm 2^{n-1}\}$. Hence, by Theorem~\ref{T.char}(b),
$\psi_n\in\M_{n-2}$.

Since $n\geq 3$, we have constructed examples of wavelets in each of the
equivalence classes $\M_n$, $n\geq 1$.
\hfill\qed

\section{A family of wavelets belonging to the class \protect{$\caM_{\bf 0}$}}
\label{sec:equiv2}

In this section we give examples of a family of wavelets in $\M_0$. In
fact, we shall show that all non-MSF wavelets characterized in
Theorem~\ref{thm:sn} belong to $\M_0$.

For $n\geq 3$, let $a_n$, $b_n$, $c_n$, $d_n$ and $e_n$ be as in
\S\ref{subsec:equiv1}. Define the function $b$ on $[e_n,b_n]$ as
follows:
\begin{equation*}
 b(\xi)=
 \begin{cases}
 {1}/{\sqrt 2} & {\rm if}~\xi\in[e_n,\pi]\\
                 0 & {\rm if}~\xi\in [\pi,b_n].
 \end{cases}
\end{equation*}
Then we extend $b$ to the whole of $S_n$ by using (i)--(iv) of
Theorem~\ref{thm:sn}. Observe that $[e_n,b_n]\cap{\rm
supp}~b\cap((1/2^{n-1})\
{\rm supp}~b)=[e_n,\pi]$. We define $\theta$ on $\R$ as
\begin{equation*}
 \theta(\xi)=
 \begin{cases}
  \pi & {\rm if}~\xi\in[e_n,\pi],\\
  0   & {\rm otherwise.}
 \end{cases}
\end{equation*}
It is easy to see that $\theta$ satisfies the functional equation in (v)
of Theorem~\ref{thm:sn}. This choice of $b$ and $\theta$ will give us
the wavelet $w_n$, where
$\widehat{w}_n(\xi)=\hbox{e}^{i\theta(\xi)}b(\xi)$. That is,
\begin{equation*}
 \widehat{w}_n(\xi)=
  \begin{cases}
  1  & {\rm if}~\xi\in[-\pi,-a_n]\cup[a_n,e_n]\cup[2^{n-1}\pi,d_n],\\[.3pc]
  {1}/{\sqrt 2} & {\rm if}~\xi\in[-d_n,-2^{n-1}\pi]
  \cup[-b_n,-\pi]\cup[c_n,2^{n-1}\pi],\\[.3pc]
  -{1}/{\sqrt 2} & {\rm if}~\xi\in[e_n,\pi],\\[.3pc]
  0 & {\rm otherwise.}
  \end{cases}
\end{equation*}

Since $[-b_n,-\pi]+2\pi=[e_n,\pi]$, and $\widehat{w}_n$ does not vanish
on the sets $[-b_n,-\pi]$ and $[e_n,\pi]$, it is clear that
$1\in{\mathcal E}(w_n)$. Hence $w_n\in\M_0$, by Theorem~\ref{T.char}(c).

\setcounter{theo}{0}
\begin{rema}
{\rm The above example is a particular case of the fact that
all wavelets characterized in Theorem~\ref{thm:sn} belong to
$\M_0\cup\M_\infty$. To see this, let $\psi$ be such a wavelet. If
$\psi$ is an MSF wavelet, then $\psi\in\M_\infty$ by the
characterization of $\M_\infty$. On the other hand, if $\psi$ is not
MSF, then there is a non-trivial set $K\subset[e_n,b_n]$ such that
$0<|\hat\psi(\xi)|<1$ for a.e. $\xi\in K$. By condition (iii) of
Theorem~\ref{thm:sn}, $0<|\hat\psi(\xi)|<1$ for a.e. $\xi\in K-2\pi$ as
well. As above, $1\in{\mathcal E}(\psi)$ proving that $\psi\in\M_0$.}
\end{rema}

\section*{Acknowledgements}

The research done by the first author was supported by the National
Board for Higher Mathematics (NBHM), Govt. of India.

\end{document}